\documentclass[reqno]{amsart}
\usepackage{mathrsfs}
\usepackage{color}
\usepackage{amsmath}
\usepackage{amsfonts}
\usepackage{amssymb}
\usepackage{graphicx}
\usepackage{hyperref}

 \newtheorem{Theorem}{Theorem}[section]

 \newtheorem{Question}[Theorem]{Question}

 \numberwithin{equation}{section}
\begin{document}

\title[A remark on the extension of $L^{2}$ holomorphic functions]
 {A remark on the extension of $L^{2}$ holomorphic functions}

\author{Qi'an Guan}
\address{Qi'an Guan: School of
Mathematical Sciences, Peking University, Beijing 100871, China.}
\email{guanqian@math.pku.edu.cn}

\thanks{The authors were partially supported by NSFC}



\date{}

\dedicatory{}

\commby{}

%%% ----------------------------------------------------------------------

\begin{abstract}
In this note, we answer a question on the extension of $L^{2}$ holomorphic functions posed by Ohsawa.
\end{abstract}

%%% ----------------------------------------------------------------------
\maketitle
%%% ----------------------------------------------------------------------

\section{an answer to a question posed by Ohsawa}
In \cite{Ohsawa2017}, Ohsawa gave a survey on a recent ''remarkable'' progress (c.f. \cite{B-L,Bl_inv,G-ZhouL2_CR,G-ZhouL2_ann,G-ZhouL2_Sci,G-Zhou-ZhuL2_CR})
around the famous Ohsawa-Takegoshi $L^{2}$ extension theorem \cite{O-T}.
After that,
Ohsawa recalled the following consequence of the main result in \cite{G-ZhouL2_ann}, 
and presented a shorter proof based on a general criterion for the extendibility in \cite{Ohsawa2017}.

\begin{Theorem}
\label{coro:GZ-domain}(\cite{G-ZhouL2_ann}, see also Theorem 0.1 in \cite{Ohsawa2017})
Let $D\subset\mathbb{C}^n$ be a pseudoconvex domain, and let $\varphi$ be a plurisubharmonic
function on $D$ and $H=\{z_n=0\}$. Then for any holomorphic
function $f$ on $H$ satisfying
$$\int_{H}|f|^{2}e^{-\varphi}dV_{H}<\infty,$$
there exists a holomorphic function $F$ on $D$ satisfying $F = f$ on $H$ and
\begin{eqnarray*}
\int_{D}|F|^{2}e^{-\varphi-(1+\varepsilon)\log(1+|z_{n}|^{2})}dV_{D}
\leq\frac{\pi}{\varepsilon}\int_{H}|f|^{2}e^{-\varphi}dV_{H}.
\end{eqnarray*}
\end{Theorem}

In \cite{Ohsawa2017}, 
considering general plurisubharmonic function $\psi(z_{n})$ instead of $(1+\varepsilon)\log(1+|z_{n}|^{2})$
in Theorem \ref{t:ohsawa2017}, 
Ohsawa posed the following question on the extension of $L^{2}$ holomorphic functions.

\begin{Question}
\label{Q:ohsawa2017}
Given a subharmonic function $\psi$ on $\mathbb{C}$ such that $\int_{\mathbb{C}}e^{-\psi}<+\infty$,
for any subharmonic function $\varphi$ on $\mathbb{C}$,
can one find a holomorphic function $f$ on $\mathbb{C}$ satisfying $f(0)=1$,
and
$$\int_{\mathbb{C}}|f|^{2}e^{-\varphi-\psi}\leq e^{-\varphi(0)}\int_{\mathbb{C}}e^{-\psi}?$$
\end{Question}

When $\psi$ does not depend on $\arg z$, Sha Yao gave a positive answer to Question \ref{Q:ohsawa2017} in her Ph.D thesis 
by using the main result in \cite{G-ZhouL2_ann}.

In the present article, 
we give the following (negative) answer to Question \ref{Q:ohsawa2017}.

\begin{Theorem}
\label{t:ohsawa2017}
There exist subharmonic functions $\psi$ and $\varphi$ on $\mathbb{C}$ satisfying

$(1)$ $\int_{\mathbb{C}}e^{-\psi}<+\infty$;

$(2)$ $\varphi(0)\in(-\infty,+\infty)$;

$(3)$ for any holomorphic function $f$ on $\mathbb{C}$ satisfying $f(0)=1$, 
$\int_{\mathbb{C}}|f|^{2}e^{-\varphi-\psi}=+\infty$ holds.
\end{Theorem}

\section{Proof of Theorem \ref{t:ohsawa2017}}
Let  $\psi=2\max\{c_{1}\log|z-1|,c_{2}\log|z-1|\}$ and $\varphi=(1-c_{1})\log|z-1|$,
where $c_{1}\in(\frac{1}{2},1)$ and $c_{2}\in(1,\frac{3}{2})$.

We prove Theorem \ref{t:ohsawa2017} by contradiction:
if not, then there exists holomorphic function $f$ on $\mathbb{C}$ satisfying $f(0)=1$,
and
\begin{equation}
\label{equ:ohsawa1}
\int_{\mathbb{C}}|f|^{2}e^{-\varphi-\psi}<+\infty.
\end{equation}

Note that $(\psi+\varphi)|_{\{|z|<1\}}=2\log|z-1|$,
then inequality \eqref{equ:ohsawa1} implies that $f(1)=0$.

Note that $\psi+\varphi-2(1-c_{1}+c_{2})\log|z|$ is bounded near $\infty$,
then inequality \eqref{equ:ohsawa1} implies that $f$ is a polynomial.
Furthermore, it follows from $1-c_{1}+c_{2}<2$ and inequality \eqref{equ:ohsawa1} that the degree of $f$ must be $0$,
which contradicts $f(1)=0$.
This proves the present theorem.

%%%------------------------------------------------------------------------

\vspace{.1in} {\em Acknowledgements}. The author would like to thank Professor Takeo Ohsawa for giving us series talks in Peking University and sharing his recent work.

\bibliographystyle{references}
\bibliography{xbib}

\begin{thebibliography}{100}
\bibitem{B-L}B. Berndtsson and L. Lempert, A proof of the Ohsawa-Takegoshi theorem with sharp estimates. J. Math. Soc. Japan 68 (2016), no. 4, 1461--1472.
\bibitem{Bl_inv}Z. Blocki, Suita conjecture and the Ohsawa-Takegoshi extension theorem, Invent.
Math. 193 (2013) 149--158.
\bibitem{G-ZhouL2_CR}Q.A. Guan, X.Y. Zhou, Optimal constant problem in the $L^2$ extension theorem. C. R. Math. Acad. Sci. Paris 350 (2012), no. 15--16, 753--756.
\bibitem{G-ZhouL2_ann}Q.A. Guan, X.Y. Zhou, A solution of an $L^2$ extension problem with an optimal estimate and applications, Ann. of Math. (2) 181 (2015), no. 3, 1139--1208.
\bibitem{G-ZhouL2_Sci}Q.A. Guan, X.Y. Zhou, Optimal constant in an $L^2$ extension problem and a proof of a conjecture of
Ohsawa, Sci. China Math. 58(1) (2015) 35--59.
\bibitem{G-Zhou-ZhuL2_CR}Q.-A. Guan, X.-Y. Zhou and L. Zhu, On the Ohsawa¨CTakegoshi L2 extension theorem
and the twisted Bochner¨CKodaira identity, C. R. Math. Acad. Sci. Paris 349(13¨C14)
(2011) 797--800.
\bibitem{Ohsawa2017}T. Ohsawa, On the extension of $L^2$ holomorphic functions VIII--a remark on a theorem of Guan and Zhou. Internat. J. Math. 28 (2017), no. 9, 1740005, 12 pp.
\bibitem{O-T} T. Ohsawa, K. Takegoshi, \emph{On the extension of $L^2$ holomorphic functions}, Math. Z. 195 (1987), 197--204.
\end{thebibliography}

\end{document}